\documentclass[12pt]{article}%
\usepackage{amsmath}
\usepackage{graphicx}%
\usepackage{amsfonts}%
\usepackage{amssymb}
\newtheorem{theorem}{Theorem}[section]

\begin{document}

\title{Obstructions to Deformations of D.G. Modules}
\author{Trina Armstrong and Ron Umble\thanks{This paper reports the results of an
undergraduate honors project directed by the second author.}}
\date{Submitted January 19, 1995\\
Revised June 23, 2001}
\maketitle

\begin{abstract}
Let $\mathbf{k}$ be a field and $n\geq1$. There exists a differential graded
$\mathbf{k}$-module $(V,d)$ and various approximations to a differential
$d+td_{1}+t^{2}d_{2}+\cdots+t^{n}d_{n}$ on $V[[t]],$ one of which is a
non-trivial polynomial deformation, another is obstructed, and another is
unobstructed at order $n$. The analogous problem in the category of
$\mathbf{k}$-algebras in characteristic zero remains a long-standing open question.

\end{abstract}

\section{Introduction}

Most deformation theories, including those of
Froelicher-Nijenhuis-Kodaira-Spencer for complex analytic manifolds and of
Gerstenhaber for algebras, introduce both a primary obstruction to extending
an infinitesimal deformation as well as ones of successively higher order,
which appear after each previous obstruction is passed. Somewhat surprisingly,
examples of infinitesimals with a non-vanishing higher order obstruction are
often difficult to find. And this difficulty seems to increase with richer
algebraic structure. In the deformation theory of Hopf algebras--a setting
with rich algebraic structure--S. D. Schack's ``primary obstruction
conjecture'' asserts that an infinitesimal extends to a deformation whenever
its primary obstruction vanishes. In the more relaxed setting of finite
dimensional associative algebras, Gerstenhaber and Schack \cite{Gersten3}
cleverly apply topological methods to produce infinitesimals with vanishing
primary obstructions on certain rigid algebras in characteristic $p>0$. But
these examples are by no means simple, and the analogous problem in
characteristic zero remains open (\cite{Gersten1}, p.61).

In this paper we consider the deformation theory of differential graded
$\mathbf{k}$-modules (\textit{d.g.m.'s}) over an arbitrary field $\mathbf{k.}$
In this setting, the algebraic structure is quite simple and infinitesimal
deformations with non-vanishing higher order obstructions are easily observed.
Given an integer $n\geq1$, we exhibit a d.g.m.\textit{\ }$(V,d)$ and various
approximations to a differential $d+td_{1}+t^{2}d_{2}+\cdots+t^{n}d_{n}$ on
$V[[t]],$ one of which is a non-trivial polynomial deformation, another is
obstructed, and another is unobstructed at order $n$.

Classically, one finds the deformation theory of d.g.m.'s imbedded in richer
theories such as the deformation theory of differential graded $\mathbf{k}%
$-algebras; as such, it has been considered by many authors (see for example
\cite{Burghelea}, \cite{Felix}, \cite{Halperin}, \cite{Markl}, \cite{Umble}).
Recently, Gerstenhaber and Wilkerson \cite{Gersten2} gave a systematic
treatment of the deformation theory of d.g.m.'s, whose relative simplicity
provides an ideal entry point for the uninitiated reader. Many of the ideas
and techniques familiar in more general settings appear here, but without the
tight constraints that typically confound even the most basic computational
examples. We begin with a review of the ideas we need.\bigskip\ 

\section{Cohomology and Deformations of D.G.\newline Modules}

Let $R$ be a commutative ring with identity---in this discussion $R$ will be
either a field $\mathbf{k}$ or the ring of formal power series $\mathbf{k}%
[[t]]=\{\sum\lambda_{i}t^{i}\mid\lambda_{i}\in\mathbf{k}\}.$ Let
$\{M_{p}\}_{p\in\mathbf{Z}}$ be a sequence of $R$-modules. Then $M=\sum
_{p\in\mathbf{Z}}M_{p}$ is a \textit{graded }$R$-module; an element $x_{p}\in
M_{p} $ is said to be \textit{homogeneous of degree} $p,$ in which case we
write $\left|  x_{p}\right|  =p$. Let $M$ and $N$ be graded $R$-modules. An
$R$-linear map $f:M\rightarrow N$ has \textit{degree} $p$ if $\left|
f(x)\right|  =\left|  x\right|  +p$ for each homogeneous $x\in M,$ in which
case we write $\left|  f\right|  =p.$\ The set of all such maps of degree $p$
is denoted by $Hom_{R}^{p}(M,N).$ A $R$-linear map $d:M\rightarrow M$ of
degree $\pm1$ such that $d^{2}=0$ is a \textit{differential} on $M$ and the
pair $(M,d)$ is a \textit{differential graded }$R$\textit{-module
}(\textit{d.g.m.})$.$

Let $(V,d_{V})$ and $(M,d_{M})$ be a d.g.m.\textit{'s}. For each
$p\in\mathbf{Z},$ define the $R$-module of $p$\textit{-cochains on }$V$
\textit{with coefficients in} $M$ by $C^{p}(V;M)=$\linebreak $Hom_{R}%
^{-p}(V,M).$ Let $C^{\ast}(V;M)=\sum_{p\in\mathbf{Z}}C^{p}(V;M)$ and define
the $p^{th}$ \textit{coboundary }map $\delta^{p}:C^{p}(V;M)\rightarrow
C^{p+1}(V;M)$ by
\[
\delta^{p}(f)=d_{M}f-\left(  -1\right)  ^{p}fd_{V},
\]
then%
\[
\delta=\sum_{p\in\mathbf{Z}}\delta^{p}\in Hom_{\mathbf{k}}^{1}\left(  C^{\ast
}\left(  V;M\right)  ,C^{\ast}\left(  V;M\right)  \right)
\]
and it trivial to check that $\delta^{2}=0$ so that $\left\{  C^{\ast}\left(
V;M\right)  ,\delta\right\}  $ is a d.g.m. The \textit{cohomology of }$V$
\textit{with coefficients in }$M$ is the cohomology of $\left\{  C^{\ast
}\left(  V;M\right)  ,\delta\right\}  ,$ i.e., the graded $R$-module $H^{\ast
}\left(  V;M\right)  =\sum_{p\in\mathbf{Z}}H^{p}\left(  V;M\right)  ,$ where
$H^{p}\left(  V;M\right)  =\ker\left(  \delta^{p}\right)  /$Im$\left(
\delta^{p-1}\right)  .$ The elements of $\ker\left(  \delta^{p}\right)  $ are
the $p$-\textit{cocycles}; the elements of Im$\left(  \delta^{p-1}\right)  $
are the $p$\textit{-coboundaries. }Two $p$-cocycles $f$ and $g$ are
\textit{cohomologous }provided $f-g\in$Im$\left(  \delta^{p-1}\right)  ;$ the
class $[f]\in H^{p}\left(  V;M\right)  $ is called the \textit{cohomology
class of} $f.$ When $V=M$ we let $C^{\ast}\left(  V\right)  =C^{\ast}\left(
V;V\right)  $ and $H^{\ast}\left(  V\right)  =H^{\ast}\left(  V;V\right)  .$

Let $\left(  V,d\right)  $ be a d.g.m. over a field $\mathbf{k.}$ The
cohomology $H^{\ast}(V)$ directs the deformation theory of $(V,d)$ in the
following way: Let $t$ be an indeterminant of degree $0$ and imbed $V$ as the
set of constants in the $\mathbf{k}[[t]]$-module of formal power series
$V[[t]]=\{\sum t^{i}v_{i}\mid v_{i}\in V\}$. Extend $d$ to $V[[t]]$ via
$d(\sum t^{i}v_{i})=\sum t^{i}d(v_{i})$---this is the unique $\mathbf{k}%
[[t]]$-linear extension of $d$---and obtain the \textit{d.g.} $\mathbf{k}%
[[t]]$-module $V_{0}[[t]]=(V[[t]],d).$ Given $\left\{  d_{i}\in C^{1}%
(V)\right\}  _{i\geq1},$ extend each $d_{i}$ to $V[[t]]$ in the same way we
extended $d$ and define
\[
d_{t}=d+td_{1}+t^{2}d_{2}+t^{3}d_{3}+\cdots.
\]
If $d_{t}$ is a differential on $V[[t]]$, then $V_{t}=(V[[t]],d_{t})$ is a
\textit{deformation} \textit{of} $(V,d)$ \textit{as a} \textit{d.g.
}$\mathbf{k}[[t]]$-\textit{module}. The \textit{trivial deformation}
$V_{t}=V_{0}[[t]]$ satisfies $d_{t}=d$, i.e., $d_{i}=0$ for all $i.$ A
deformation $V_{t}$ with differential $d_{t}$ such that $d_{n}\neq0$ and
$d_{i}=0$ for all $i>n$ is called a \textit{polynomial deformation of order\ }$n.$

Let $V_{t}=(V[[t]],d_{t})$ be a deformation. Define $\mathcal{O}_{0}=0$ and
for each $n\geq1$ consider the $2$-cochain
\[
\mathcal{O}_{n}=-\sum_{i=1}^{n}d_{i}d_{n-i+1}.
\]
Expand the right hand side of $d_{t}\circ d_{t}=0$ and equate coefficients to
obtain the relations
\begin{equation}
\{\delta(d_{n+1})=\mathcal{O}_{n}\}_{n\geq0}. \label{relations}%
\end{equation}
Observe that $d_{1}$ is a $1$-cocycle and that $\mathcal{O}_{n}$ is a
cobounding $2$-cocycle for each $n;$ these are necessary and sufficient
conditions for an arbitrary series $d_{t}^{\prime}=d+td_{1}^{\prime}%
+t^{2}d_{2}^{\prime}+\cdots$ to be a differential, i.e., for $(V[[t]],d_{t}%
^{\prime})$ to be a deformation. The differential $d_{t}$ is commonly referred
to as a deformation of $d;$ the $1$-cocycle $d_{1}$ is commonly referred to as
an \textit{infinitesimal }deformation. The $2$-cocycles $\{\mathcal{O}%
_{n}\}_{n\geq1}$ are the \textit{obstructions }to extending the linear
approximation $d+td_{1}$; and $\mathcal{O}_{1}$ is called the \textit{primary} obstruction.

This suggests the following inductive strategy for constructing deformations:
Given a d.g.m.\textit{\ }$(V,d),$ choose a $1$-cocycle $d_{1}$ as required by
the $0^{th}$ relation in (\ref{relations}) and obtain the linear approximation
$d+td_{1}$. Inductively, assume that for each $k\leq n,$ some $1$-cochain
$d_{k}$ has been chosen so that $\delta(d_{k})=\mathcal{O}_{k-1}.$ If
$\mathcal{O}_{n}$ fails to cobound, the approximation $d+td_{1}+\cdots
+t^{n}d_{n}$ \textit{is obstructed at order }$n$ and the process terminates.
Otherwise, choose a $1$-cochain $d_{n+1}$ such that $\delta(d_{n+1}%
)=\mathcal{O}_{n}$ and extend the approximation to $d+td_{1}+\cdots
+t^{n+1}d_{n+1}.$ If the process can be continued indefinitely, there is a
differential $d_{t}=d+td_{1}+t^{2}d_{2}+\cdots$ whose partial sums are the
approximations obtained inductively. Indeed the inductive process can be
continued indefinitely if $\mathcal{O}_{n}\in\lbrack0]$ in $H^{2}(V)$ for all
$n.$ Since this happens automatically whenever $H^{2}(V)=0$ we obtain:\medskip

\begin{theorem}
Let $(V,d)$ be a d.g.m.\textit{\ }and let $d_{1}$ be a $1$-cocycle. If
$H^{2}(V)=0,$ there exists a sequence of $1$-cochains $\{d_{i}\}_{i\geq2}$
such that $(V[[t]],d_{t}=d+td_{1}+t^{2}d_{2}+\cdots)$ is a deformation.\medskip
\end{theorem}

Two deformations $V_{t}=(V[[t]],d_{t})$ and $V_{t}^{\prime}=(V[[t]],d_{t}%
^{\prime})$ are \textit{equivalent} if there exists a $\mathbf{k}[[t]]$-linear
automorphism $\phi_{t}:V[[t]]\rightarrow V[[t]]$ such that

\begin{enumerate}
\item $\phi_{0}=Id_{V[[t]]}$ and

\item $d_{t}\phi_{t}=\phi_{t}d_{t}^{\prime}.$
\end{enumerate}

\noindent Condition (1) implies the existence of maps $\{\phi_{i}\in
C^{0}(V)\}_{i\geq1}$ such that $\phi_{t}=Id+t\phi_{1}+t^{2}\phi_{2}+\cdots.$
Condition (2), the naturality condition, implies that $\phi_{t}$ induces an
$\mathbf{k}[[t]]$-linear isomorphism of d.g.m.'s $H^{\ast}(V[[t]],d_{t}%
^{\prime})\approx H^{\ast}(V[[t]],d_{t}).$ When $\phi_{t}$ exists we call it
an \textit{equivalence} and write $\phi_{t}:V_{t}\sim V_{t}^{\prime}.$

Suppose that $V_{t}$ and $V_{t}^{\prime}$ are deformations of $(V,d)$ and
assume that $\phi_{t}:V_{t}\sim V_{t}^{\prime}$ is an equivalence. Expand and
collect first order coefficients in relation (2) and obtain
\[
d_{1}^{\prime}-d_{1}=d\phi_{1}-\phi_{1}d=\delta(\phi_{1});
\]
thus $d_{1}$ and $d_{1}^{\prime}$ are cohomologous. In particular, if
$V_{t}^{\prime}$ is the trivial deformation $V_{0}[[t]],$ then $d_{1}^{\prime
}=0$ in which case
\begin{equation}
\delta(\phi_{1})=-d_{1}. \label{coboundary}%
\end{equation}
A d.g.m.\textit{\ }$(V,d)$ is \textit{rigid} if every deformation $V_{t}$ is
equivalent to the trivial deformation $V_{0}[[t]].$ We conclude this section
with a standard theorem whose proof is included for completeness.\medskip

\begin{theorem}
Let $(V,d)$ be d.g.m.\textit{\ }If $H^{1}(V)=0,$ then $(V,d)$ is rigid.\medskip
\end{theorem}

\textit{Proof}: Let $V_{t}=(V[[t]],d_{t}=d+td_{1}+t^{2}d_{2}+\cdots)$ be a
deformation; by assumption, there exists a $0$-cochain $\phi_{1}$ such that
$\delta(\phi_{1})=-d_{1}.$ Consider the $\mathbf{k}[[t]]$-linear isomorphism
$\phi_{t}^{(1)}=Id-t\phi_{1}$ and define $d_{t}^{(1)}=\phi_{t}^{(1)}d_{t}%
[\phi_{t}^{(1)}]^{-1}.$ Then $\phi_{t}^{(1)}:(V[[t]],d_{t}^{(1)})\sim V_{t}$
is an equivalence. Expanding $d_{t}^{(1)}\phi_{t}^{(1)}=$ $\phi_{t}^{(1)}%
d_{t}$ and equating first order coefficients gives $d_{1}^{(1)}=(d\phi
_{1}-\phi_{1}d)+d_{1}=0,$ so the linear term in $d_{t}^{(1)}$ vanishes.
Inductively, suppose that $V_{t}^{(r)}=(V[[t]],d_{t}^{(r)})$ is a deformation
with $d_{i}^{(r)}=0$ for all $i\leq r$ and that $\phi_{t}^{(r)}=\Pi_{i=1}%
^{r}(Id-t^{i}\phi_{i}):V_{t}^{(r)}\sim V_{t}$ is an equivalence. Then
$\delta(d_{r+1}^{(r)})=0,$ so by assumption there exists a $0$-cochain
$\phi_{r+1}$ such that $\delta(\phi_{r+1})=-d_{r+1}^{(r)}.$ Define $\phi
_{t}^{(r+1)}=\Pi_{i=1}^{r+1}(Id-t^{i}\phi_{i})$ and $d_{t}^{(r+1)}=\phi
_{t}^{(r+1)}d_{t}^{(r)}[\phi_{t}^{(r+1)}]^{-1}$; then $d_{r+1}^{\left(
r+1\right)  }=\left(  d\phi_{r+1}-\phi_{r+1}d\right)  +d_{r+1}^{\left(
r\right)  }=0$ so that $d_{i}^{(r+1)}=0$ for all $i\leq r+1$ and $\phi
_{t}^{(r+1)}:(V[[t]],d_{t}^{(r+1)})\sim V_{t}$ is an equivalence$.$ Hence
there is a sequence of equivalences $\{\phi_{t}^{(r)}=\Pi_{i=1}^{r}%
(Id-t^{i}\phi_{i})\}_{r\geq1}$ that $t$-adically converge to an equivalence
$\phi_{t}^{(\infty)}=\Pi_{r\geq1}(Id-t^{r}\phi_{r}):V_{0}[[t]]\sim
V_{t}.\mathit{\ }$\bigskip\ 

Analogs of Theorem 1 and Theorem 2 appear in the deformation theory of
differential graded algebras with an apparent shift of dimension. This
dimension shift, which reflects nothing more than a change in point-of-view,
is discussed in \cite{Umble}. We are ready for the construction promised.

\section{Examples of Deformations and Obstructed Approximations}

Let $\mathbf{k}$ be any field, let $V=\sum_{p\in\mathbf{Z}}V_{p}$ be a graded
$\mathbf{k}$-module and choose a basis $\{x_{\alpha}\}$ for $V.$ If $\left\{
x_{\beta}\right\}  $ is a set of vectors in $V$, let $\left\langle x_{\beta
}\right\rangle $ denote the $\mathbf{k}$-linear span. Consider the associated
graded $\mathbf{k}$-module $Hom_{\mathbf{k}}^{\ast}\left(  V,V\right)
=\left\langle x_{i}\frac{\partial}{\partial x_{j}}\right\rangle ,$ where
\[
x_{i}\frac{\partial x_{k}}{\partial x_{j}}=\left\{
\begin{array}
[c]{cl}%
x_{i}, & j=k\\
0 & \text{otherwise.}%
\end{array}
\right.
\]
Thus
\[
x_{i}\frac{\partial}{\partial x_{j}}\left(  x_{k}\frac{\partial}{\partial
x_{\ell}}\right)  =x_{i}\frac{\partial x_{k}}{\partial x_{j}}\frac{\partial
}{\partial x_{\ell}}.
\]

In particular, let $V_{p}=\{0\}$ for all $p\leq0$ and $V_{p}=\left\langle
x_{2p-1},x_{2p}\right\rangle $ for all $p\geq1;$ define%
\[
d=\sum_{i=1}^{\infty}x_{6i-5}\frac{\partial}{\partial x_{6i-3}}.
\]
Clearly $d^{2}=0$ so that $(V,d)$ is a \textit{d.g.m}. Fix a positive integer
$n\geq2$ and define
\[
d_{1}=x_{1}\frac{\partial}{\partial x_{4}}+\sum_{i=1}^{n-1}x_{6i-2}%
\frac{\partial}{\partial x_{6i}}.
\]
Since $\delta(d_{1})=0,$ the expression $d+td_{1}$ is a linear approximation
to a deformation. Furthermore,
\[
\delta(-x_{3}\frac{\partial}{\partial x_{6}})=-x_{1}\frac{\partial}{\partial
x_{6}}=-d_{1}^{2}=\mathcal{O}_{1}%
\]
so the primary obstruction vanishes in cohomology. For $2\leq k\leq n$ define
\begin{equation}
d_{k}=-x_{6k-9}\frac{\partial}{\partial x_{6k-6}}+(1-\delta_{n,k}%
)x_{6k-5}\frac{\partial}{\partial x_{6k-2}}, \label{dk}%
\end{equation}
where $\delta_{n,k}$ is the Kronecker delta, and consider the subspaces
$S=\left\langle x_{2i}\right\rangle _{i\geq1}$ and $S^{\bot}=\left\langle
x_{2i-1}\right\rangle _{i\geq1}$ spanned by basis elements of even and odd
index, respectively. Note that $d_{i}$ is supported on $S$ for $1\leq i\leq
n;$ and furthermore, $d_{i}$ takes values in $S^{\bot}$ for $2\leq i\leq n;$
hence $d_{i}d_{j}=0$ for $2\leq i,j\leq n$ and $d_{1}d_{k}=0$ for $2\leq k\leq
n.$ Therefore
\[
\mathcal{O}_{k}=-d_{k}d_{1}=\left\{
\begin{array}
[c]{cc}%
-x_{6k-5}\frac{{\LARGE \partial}}{{\Large \partial x}_{{\large 6k}}}%
\medskip, & 2\leq k<n\\
0, & n=k.
\end{array}
\right.
\]
On the other hand, $d_{k}d=0$ since $d(S)\subset S^{\bot},$ hence%
\[
\delta(d_{k})=dd_{k}.
\]
It is now a simple matter to check that
\[
\delta(-d_{k+1})=-x_{6k-5}\frac{\partial}{\partial x_{6k}}=\mathcal{O}_{k}%
\]
for $2\leq k<n.$ Now for all $i>n,$ set $d_{i}=0$ so that $\mathcal{O}%
_{i}=-d_{i}d_{1}=0$ and obtain a polynomial deformation $V_{t}%
=(V[[t]],\ d+td_{1}+t^{2}d_{2}+\cdots+t^{n}d_{n})$ of $(V,d)$. To establish
non-triviality, assume that $\phi_{t}:V_{t}\sim V_{0}[[t]].$ Then there exists
a cochain $\phi_{1}\in C^{0}(V)$ such that $\delta\left(  \phi_{1}\right)
=-d_{1}$, i.e.,
\[
d\phi_{1}-\phi_{1}d=-x_{1}\frac{\partial}{\partial x_{4}}-x_{4}\frac{\partial
}{\partial x_{6}}+\text{\textit{(other terms)}}.
\]
Now $d\left(  x_{6}\right)  =0$ so $\delta\left(  \phi_{1}\right)  \left(
x_{6}\right)  =d\phi_{1}\left(  x_{6}\right)  =-x_{4}.$ But this is a
contradiction since $\phi_{1}\left(  x_{6}\right)  \in V_{3}=\left\langle
x_{5},x_{6}\right\rangle $ and $d\left(  V_{3}\right)  =0.$ Finally, if we
redefine $d_{1}=x_{4}\frac{\partial}{\partial x_{6}}$ then $d_{t}=d+td_{1}$ is
a non-trivial linear deformation, as the reader can easily check. At the other
extreme, redefine $d_{1}=x_{1}\frac{\partial}{\partial x_{4}}+\sum
_{i=1}^{\infty}x_{6i-2}\frac{\partial}{\partial x_{6i}}$ and obtain a
non-trivial deformation with non-zero terms of all orders. We have proved:\medskip

\begin{theorem}
For each $n\geq1,$ there exists a d.g.m $(V,d)$ and a non-trivial polynomial
deformation of order $n.$ Furthermore, a sequence $\left\{  d_{i}\right\}
_{i\geq1}$ can be chosen such that $V_{t}=\left(  V[[t]],\ d_{t}%
=d+td_{1}+t^{2}d_{2}+\cdots\right)  $ is a non-trivial deformation with
$d_{i}\neq0$ for all $i.$\medskip
\end{theorem}

We conclude by constructing an approximation to a differential on $V[[t]]$
that is obstructed at order $n.$ First assume that $2\leq k<n$ and define
$d_{1}$ and $d_{k}$ as in the example above. Extend $d+td_{1}$ to an
approximation of order $n-1$ in the same manner as before and obtain
$\mathcal{O}_{n-1}=-d_{n-1}d_{1}=-x_{6n-11}\frac{\partial}{\partial x_{6n-6}%
}.$ Define
\[
d_{n}=-x_{6n-9}\frac{\partial}{\partial x_{6n-6}}+x_{6n-6}\frac{\partial
}{\partial x_{6n-4}}
\]
and observe that $\delta(d_{n})=\mathcal{O}_{n-1}$ so that the approximation
of order $n-1$ extends to order $n.$ This time,
\[
\mathcal{O}_{n}=-d_{1}d_{n}=-x_{6n-8}\frac{\partial}{\partial x_{6n-4}}
\]
so that
\[
\mathcal{O}_{n}\left(  V\right)  \subset V_{3n-4}.
\]
Since $d(V)\subset\sum_{i=1}^{\infty}V_{3p-2},$ we have
\[
\delta(f)\left(  V\right)  =\left(  df+fd\right)  \left(  V\right)
\subset\sum_{p=1}^{\infty}V_{3p-2}\oplus V_{3p},
\]
for all $f\in C^{1}\left(  V\right)  .$ Consequently $\mathcal{O}_{n}$ fails
to cobound and the approximation is obstructed at order $n.$ Finally, when
$n=1,$ define $d_{1}=x_{4}\frac{\partial}{\partial x_{6}}+x_{6}\frac{\partial
}{\partial x_{8}}$ and observe that $\delta(d_{1})=0$ while $\mathcal{O}%
_{1}=-d_{1}^{2}=-x_{4}\frac{\partial}{\partial x_{8}}$ fails to cobound since
$\mathcal{O}_{1}\left(  V\right)  \subset V_{2}$. This proves our main result:\medskip

\begin{theorem}
For each $n\geq1,$ there exists a differential graded $\mathbf{k}$-module
$(V,d)$ and an obstructed approximation to a differential $d+td_{1}+t^{2}%
d_{2}+\cdots t^{n}d_{n}$ on $V[[t]].$\bigskip
\end{theorem}

\noindent{\LARGE Acknowledgement}\bigskip

We wish to thank Murray Gerstenhaber for his helpful comments and suggestions.

\bigskip

\noindent Department of Health Evaluation Services, Penn State college of
Medicine, MC H173, PO\ Box 850, 500 University Dr., Hershey, PA 17033.
Phone:\ (717) 531-7178; e-mail: tja3@psu.edu\bigskip

\noindent Department of Mathematics, Millersville University of Pennsylvania,
Millers- ville, PA 17551. Phone: (717) 872-3708; e-mail: Ron.Umble@millersville.edu\bigskip

\noindent Manuscript number:\ PROC6293
\end{document}